\documentclass[11pt]{article}
\usepackage{epsfig,amssymb}
\usepackage{verbatim}
\usepackage{nameref}
\usepackage{hyperref}
\usepackage{sectsty}
\usepackage{enumerate}
\usepackage{epstopdf}

\newcommand{\sta}{{\rm sta}}

\newcommand{ \PP}{{P}}

\newcommand{\vsig}{\varsigma}

\newcommand{\barvsig}{\bar{\varsigma}}

\newcommand{\DD}{D}

\newcommand{\ww}{w}

\newcommand{\TT}{T}

\newcommand{\real}{{\mathbb R}} %\newcommand{\real}{{\bf R}}
\newcommand{\half}{\frac{1}{2}}

\newcommand{\bsig}{\mbox{\boldmath$\sigma$}}
\newcommand{\bepsilon}{\mbox{\boldmath$\epsilon$}}
\newcommand{\bvsig}{{\mbox{\boldmath$\varsigma$}}}

\newcommand{\barbvsig}{\bar{\mbox{\boldmath$\varsigma$}}}

\newcommand{\bxi}{{\mbox{\boldmath$\xi$}}}

\newcommand{\la}{\langle}
\newcommand{\ra}{\rangle}

\newcommand{\eba}{\begin{array}}
\newcommand{\eea}{\end{array}}
\newcommand{\ebe}{\begin{eqnarray}}%[section]
\newcommand{\eee}{\end{eqnarray}}%[section]
\newcommand{\eb}{\begin{equation}}%[section]
\newcommand{\ee}{\end{equation}}%[section]

\newcommand{\calW}{{\cal{W}}}
\newcommand{\calP}{{\cal{P}}}

\newcommand{\WW}{{W}}

\newcommand{\bG}{{\bf G}}
\newcommand{\bH}{{\bf H}}

\newcommand{\bff}{{\bf f}}

\newcommand{\bR}{{\bf R}}
\newcommand{\bv}{{\bf v}}
\newcommand{\bw}{{\bf w}}

\newcommand{\bx}{{\bf x}}

\newcommand{\bI}{{\bf I}}

\newcommand{\calS}{{\cal S}}

\newcommand{\calC}{{\cal C}}

\newcommand{\calE}{{\cal E}}

\newcommand{\calX}{{\cal X}}

\newcommand{\barbx}{\bar{\bf x}}

\newcommand{\dt}{{\,\mbox{d}t}}

\newcommand{\dx}{\mbox{ d}x}

\newcommand{\alp}{{\alpha}}

\newcommand{ \Lam}{{\Lambda}}

\newcommand{ \lam}{{\lambda}}
\newcommand{ \xx}{{ x}}

\newtheorem{thm}{Theorem}
\newtheorem{definition}{Definition}

\newcommand\FF{{F}}

\renewcommand\barbx{{\bar{\bx}}}

\newcommand\G{{G\^{a}teaux} }

\renewcommand\eb{\begin{equation}}
\renewcommand\ee{\end{equation}}

\newcommand{\VV}{V}

\newtheorem{lemma}{Lemma}
\def\bbR{\mathbb{R}}
\def\beq{\begin{equation}}
\def\eeq{\end{equation}}
\def\btheo{\begin{theorem}}
\def\etheo{\end{theorem}}
\newcommand{\mcal}{\mathcal}
\newtheorem{exam}{Example}

\def\bpmat{\begin{pmatrix}}
\def\epmat{\end{pmatrix}}

\def\bdefi{\begin{definition}}
\def\edefi{\end{definition}}

\begin{document}

\setcounter{page}{1}
$\;\;$\\
\begin{center}
{\Large \textbf{On Modelling and Complete Solutions to  General Fixpoint Problems in Multi-Scale Systems  with Applications
 }
\vspace{0.4cm}\\[0pt]}
{\textbf{Ning Ruan and David Yang Gao}}\footnote{Corresponding author. Email address:   d.gao@federation.edu.au }
\vspace{.3cm} \\
{\small {\it  Faculty of Science and Technology,
Federation University Australia,\\
 Ballarat, VIC 3353, Australia }}
$\;$\\
\end{center}

\begin{abstract}

  This paper revisits the well-studied fixed point problem from a unified  viewpoint of mathematical modeling and canonical duality theory,
  i.e. the original problem is first reformulated as a nonconvex optimization  problem,
  its  well-posedness is discussed based on objectivity principle in continuum physics;
  then the canonical duality theory is applied for solving this problem to obtain not only all   fixed points, but also their stability properties.  Applications are illustrated  by  challenging problems governed by nonconvex polynomial, exponential, and logarithmic operators.
This paper shows that within the framework of the canonical duality theory, there is no difference between   the fixed point problems
and nonconvex analysis/optimization  in multidisciplinary studies.
\end{abstract}

\section{Mathematical Modeling and Objectivity }
The fixed point problem  is a well-established subject in the area of nonlinear analysis \cite{bierlaire-Crittin2006,border1985,eaves1972}, which is usually formulated in the following
form:
\eb\label{primal0}
(\calP_0): \;\;\;\;  \bx = \FF(\bx),
\ee
where $\FF:\calX_a \rightarrow \calX_a$ is nonlinear mapping and $  \calX_a  $ is a subset of a   normed space $\calX $.
 Problem $(\calP_0)$ appears extensively in engineering  and sciences, such as   equilibrium problems, mathematical economics, game theory, and numerical methods for nonlinear dynamical systems.
   The general form of equilibrium problem was first considered by Nikaido and Isoda in 1955 as an auxiliary problem to
establish existence results for Nash equilibrium points in non-cooperative
games \cite{scarf1967,scarf-hansen1973,shellman-sikorski2002,shellman-sikorski2003}.
Mathematically speaking, the nonlinear operator  $F(\bx)$ could be any arbitrarily given vector-valued function.
 Therefore, the   fixed point problem $(\calP_0)$ is artificial. Although it can be used to ``model''  a large class of mathematical problems,
one must pay  a price:  it is impossible to develop a unified mathematical theory with powerful real-world applications.
This dilemma is due to a  gap between mathematical analysis and mathematical physics.
 Traditional methods for solving this nonlinear problem are based on linear iteration
  \cite{hirsch-papadimitriou-vavasis1989,huang-khachiyan-sikorski}.
 This paper will provide a different approach. For simplicity's sake,  we  assume that $\calX_a$ is a convex open set in $\real^n$ with
 a norm $\|\bx\|$ induced by the  bilinear form
 $\la * , * \ra :\calX \times \calX \rightarrow \real$.
\begin{lemma}
If $\FF$ is a potential operator, i.e. there exists a real-valued function $\PP:\calX_a \rightarrow \real$
 such that $\FF(\bx)  = \nabla \PP(\bx)$, then
  $(\calP_0)$ is equivalent to the following stationary point  problem:
\begin{equation}\label{primal1}
 \barbx = \arg \sta \left \{ \Pi(\bx)= \PP(\bx)-\frac{1}{2}\|\bx \|^2  \; | \;\; \forall  \bx\in \calX_a  \right \}.
\end{equation}
Otherwise,  $(\calP_0)$ is equivalent to the following global minimization problem:
\begin{equation}\label{primal2}
  \barbx = \arg \min \left \{ \Pi(\bx)= \half \| \FF(\bx) - \bx\|^2   \; | \;\; \forall  \bx \in \calX_a  \right \}.
\end{equation}
\end{lemma}
\noindent {\em Proof}.
First we assume that  $\FF(\bx)$ is potential operator, then $\bx$ is a stationary point  of $\Pi(\bx)$ if and only if $\nabla \Pi(\bx) = \nabla \PP(\bx) - \bx = 0$,
thus,  $\bx$ is also a solution to $(\calP_0)$ since $\FF(\bx) = \nabla \PP(\bx)$.

Now we assume that $\FF(\bx)$ is not a potential operator. By the fact that $\Pi(\bx) = \half \| \FF(\bx) - \bx \|^2 \ge 0 \;\;\forall \bx \in \calX$,  the vector  $\barbx$ is a global minimizer of $\Pi(\bx)$ if and only if  $\FF(\barbx) - \barbx = 0 $. Thus, $\barbx$ must be a solution to $(\calP_0)$. \hfill $\Box$\\

  By the facts that the global minimizer of an unconstrained optimization problem must be a
  stationary point, and
\eb
 \frac{1}{2} \| \FF(\bx)-\bx\|^2  = \PP(\bx) - \frac{1}{2} \|\bx\|^2,
 \;\;   \PP(\bx)=\frac{1}{2} \la  \FF(\bx) ,  \FF(\bx) \ra  -  \la \bx , \FF(\bx) \ra  +\|\bx\|^2 ,\label{eq-fp}
\ee
  the global minimization problem  (\ref{primal2})  is a special case of
  the stationary point problem (\ref{primal1}).
  Mathematically speaking, if a fixed point problem has a trivial solution, then $\FF(\bx)$ must be
  a homogeneous operator, i.e. $\FF(0) = 0$.
  For general problems,   $\FF(\bx)$ should have a nonhomogeneous term $\bff \in \real^n$. Thus,
 we can let
\begin{eqnarray}\label{pp}
\PP(\bx)=\WW(\DD \bx)-  \la \bx ,  \bff \ra ,
\end{eqnarray}
where %$\bff \in \real^n$ is a given vector,
 $\DD: \calX  \rightarrow \calW \subset \real^m$ is a linear operator,
 $\WW:\calW \rightarrow \real$ is a so-called {\em objective function}.
Objectivity is  a basic concept  in continuum physics \cite{ciarlet,holz,marsd-hugh} and
 mathematical modeling \cite{gao-aip,gao-opl}.
Its mathematical definition is  given in Gao's book (Definition 6.1.2 \cite{gao-book00}).
 \begin{definition}[Objectivity]
 Let
 $ {\cal R} $  be
 a proper orthogonal group, i.e. $\bR \in {\cal R} $ if and only if
 $  \bR^T = \bR^{-1} , \; \det \bR = 1$.
 A set $\calW_a   $ is said to be objective if
 \[
 \bR \bw \in \calW_a \;\; \forall \bw \in \calW_a, \; \forall \bR \in {\cal R}.
 \]
 A real-valued function $\WW:\calW_a \rightarrow \real$ is said to be objective if
 \eb
 \WW(\bR \bw ) = \WW(\bw) \;\; \forall \bw \in \calW_a, \; \forall \bR \in {\cal R}.
 \ee
 \end{definition}

Geometrically speaking, an objective function does not depend on rigid rotation of the system considered, but only  on certain measure of its variable.
In Euclidean space $\calW  \subset \real^m$, the simplest objective function is the $\ell_2$-norm
$\|\bw\|$ in $\real^m$  as we have
$\|{\bf R} \bw \|^2 = \bw^T {\bf R}^T {\bf R} \bw = \|\bw\|^2 \;\;  \forall {\bf R} \in {\cal R}$.
 For general $F(\bx)$, we can see  from (\ref{eq-fp}) that
$\frac{1}{2} \|\FF(\bx)\|^T  $  and $ \frac{1}{2} \|\bx\|^2$ are  objective functions. By  the fact that $\bx= F(\bx)$,
we know that  $\la \bx ,  F(\bx) \ra $ is also an
  objective function. Therefore, for a given fixed point problem, the corresponding  $\Pi(\bx)$ is naturally an objective function.

 Physically, an objective function is governed by the intrinsic  property of the system, which doesn't depend on observers. Because of Noether's theorem, the objective function $\WW(\bw)$  should be a SO($n$)-invariant and  this invariant  is equivalent to  certain  conservation law (see Section 6.1.2 \cite{gao-book00})
   Therefore, the objectivity is essential for any real-world mathematical models. It was emphasized by P. Ciarlet
that the objectivity is not an assumption, but an axiom \cite{ciarlet}.

From the viewpoint of systems theory,  if $\bx$ represents the output (or the state,  configuration, etc.) of the system, then the nonhomogeneous term $\bff$ can be viewed as the
  input (or the control, applied force, etc.), which depends on each given problem.
  %is necessary for a well-proposed problem to have nontrivial solution.
Correspondingly, the linear term $\la \bx, \bff\ra $ in (\ref{pp}) can be called the {\em subjective function} \cite{gao-aip, gao-opl}.
Thus, the fixed point problem $(\calP_0)$ can be reformulated in the following stationary point problem:
\begin{equation}\label{primal3}
(\mathcal{P} ):~~~~ \barbx = \arg \sta \left\{ \Pi(\bx)= \WW(\DD \bx)-\frac{1}{2}\|\bx\|^2 - \la \bx ,  \bff \ra   \; |
 \;\; \forall  \bx \in \calX_a  \right\}.
\end{equation}

From the  theory of nonconvex analysis, any nonconvex function can be written as a d.c. (deference of convex)
 functions \cite{jin-gao-amc}. Therefore, the fixed point problem is actually equivalent to a
 d.c. programming problem.
By the fact that $\calX$ and $\calW$ are two different spaces with different scales (dimensions),
 the  problem $(\calP)$  can be used to study general  problems in multi-scale complex systems.

For potential operator,  a fixed point is just  a stationary point, which can be easily find  by traditional linear  iteration methods.
 For nonpotential operator, the fixed point must be a global minimizer. Due to the lack of global optimality condition in traditional theory of nonlinear optimization,    to solve a general nonconvex minimization problem  is considered to be NP-hard in global optimization and computer science.
However, this paper will show that many of these nonconvex   problems can be solved in an elegant way.

\section{Properly Posted Problem and Challenges}
According to  the Brouwer fixed-point theorem we know that any continuous function from the closed unit ball in n-dimensional Euclidean space to itself must have a fixed point.
Generally speaking, for any given non-trivial input, a well-defined
system should have at least one non-trivial response.
\begin{definition}[Properly- and Well-Posed Problems]
The  problem $(\calP)$  is called
properly posed if for any given non-trivial input  $\bff \neq 0 $, it has at least one non-trivial solution. It is called
well-posed
if the solution is unique.
\end{definition}
Clearly, this definition is more general than Hadamard's well-posed problems in dynamical
systems since the continuity condition for the solution is not required.  Physically speaking, any real-world
problems should be well-posed since all natural phenomena exist uniquely. But practically, it
is difficult to model a real-world problem precisely. Therefore, properly posed problems are
allowed for the canonical duality theory. This definition is
important for understanding  challenging problems in complex systems.

 \begin{exam}[Manufacturing/Production Systems]
  {\em
  In management science, the output  is a vector $\bx \in \real^n$, which could  represent  the products of a manufacture company.
  The input  $\bff \in \real^n$ can be considered as market price (or demanding). Therefore, the subjective function
  $ \la \bx , \bff\ra = \bx^T \bff$ in this example
is  the total income of the company.
The products are produced by workers $\bw \in \real^m$. Due to the cooperation, we have $\bw = \DD \bx$ and
$\DD \in \real^{m\times n}$ is a matrix. Workers are paid by salary $\bsig = \partial \WW(\bw)$,
therefore, the objective function $\WW(\bw)$ is the cost (in this example $\WW$
 is not necessarily to be
objective since the company is a man-made system).
Let $\half \alp \|\bx\|^2$ be the profit that the company must to make, where $\alp> 0$ is a parameter,
then  $\Pi(\bx) = \WW(\DD \bx) + \half \alp \|\bx \|^2 - \bx^T \bff $ is the {\em target}  and the minimization problem $\min \Pi(\bx)$ leads to
the equilibrium equation
\[
\alp \bx  =  \bff  - \DD^T \partial_{\bw} \WW(\DD \bx) .
\]
This is a fixed point problem.
 The cost function  $\WW(\bw)$ could be  convex for a  small company, but
usually nonconvex for big companies to allow some people having  the same salaries.
}
  \end{exam}

  \begin{exam}[Lagrange Mechanics] \label{exam-lag}
  {\em
  In analytical mechanics, the configuration $\bx \in \calX  \subset \calC^1[I; \real^{n}]$ is a continuous vector-valued
  function of time $t\in I \subset \real$. Its components $\{ \xx_i \} \; (i = 1, \dots, n) $ are  known as
  the {\em Lagrangian coordinates}.  The input  $\bff(t)$   is a given  force vector function in $\real^n$. Therefore, the subjective functional in this case is
   $ \la \bx, \bff \ra = \int_I \bx(t) \cdot \bff(t) \dt$.
 The total action of the system is
  \[
 \int_I L(\bx, \dot{ \bx} ) \dt , \;\; L=  \TT(\dot{ \bx}  ) - \VV( \bx)
  \]
  where $\TT $ is the kinetic energy density,  $\VV $ is the potential density, and $L= \TT - \VV$ is the standard
  {\em Lagrangian density}. In this case, the linear operator $D = \partial_t$ is a derivative with time.
   Together, $\Pi(\bx) =\int_I [ \TT(\dot{ \bx} ) - \VV(\bx) - \bx^T \bff ] \dt   $   is called  the {\em total action}.
   For Newton mechanics, the kinetic energy   is a quadratic (objective) function
   $\TT(\bv) = \half   m \|\bv \|^2 $.
  Its  stationary condition leads to the {\em  Euler-Lagrange equation}:
   \eb
 -  m  \ddot{\bx }=  \bff + \nabla  \VV( \bx)  . \label{eq-e-l}
   \ee
 Finite difference method for solving this second-order differential equation leads to a fixed point problem \cite{ruan-gao-ima}.
 It is well-known that if the potential energy $\VV(\bx)$ is convex, the operator $\FF =  \bff + \nabla  \VV( \bx) $ is monotone and  the problem $(\calP_0)$ has stable fixed point solution.
 Correspondingly, the system has stable trajectory.
 Otherwise, the system could have chaotic solutions.
 The relation between chaos in nonlinear dynamical systems and NP-hardness in computer science  is discovered recently \cite{lat-gao-chaos}.
   } \end{exam}

   \begin{exam}[Post-Buckling of Nonlinear Gao Beam] \label{exam-gao}
  {\em
In large deformation solid mechanics, the correct nonlinear beam theory  that can be used to model post-buckling phenomenon was proposed by Gao in 1996 \cite{gao-mrc96}, which  is governed by a forth-order nonlinear differential equation:
\eb
 \chi_{xxxx} -    \frac{3}{2} \alp \chi_x^2 \chi_{xx} +  \lam \chi_{xx} =q ,\label{eq-beam}
\ee
where   $\chi(x)$ is  the deflection of the beam, which is a scaler-valued function over its domain
 $[0,L]$, where $L$ is the beam length, $\alp> 0$ is a material constant,
  the parameter  $\lam$ depends on the axial force and $q(x)$ is a given distributed lateral load.
 Clearly, this nonlinear deferential equation can be written in the following fixed point problem:
  \eb
  \chi(x) = F(x, \chi(x)), \;\; F(x,\chi) = \int_0^x \int_0^t \int_0^s \alp \left( \half \chi_x^3 - \lam \chi_x \right)ds dt dx + f(x)
  \ee
  where the function $f(x)$ depends on both the lateral load $q(x)$ and boundary conditions.
   In this case, $F(x,\chi(x))$ is a nonlinear integration operator.
This fixed point problem is equivalent to  the stationary point problem
\eb
\chi = \arg \sta \left\{ \Pi(\chi) = \int_0^L \left[ \half \chi^2_{xx} + \half \alp  \left( \half \chi^2_x -\lam\right)^2 - q(x) \right]\dx |\;\;\chi \in \calX_a \right\}. \label{eq-beamfix}
\ee
It was indicated in \cite{gao-ijnm00}   that  if $\lam< \lam_c$, the Euler buckling load  defined by
\[
 \lam_c  = \inf \frac{\int_0^L \chi^2_{xx} \dx }{\alp \int_0^L \chi^2_x \dx },
\]
the total potential $\Pi(\chi)$ is a convex functional and  the problem (\ref{eq-beamfix}) has only one fixed point. In this case,
 the beam is in pre-buckling state.
If $\lam > \lam_c$, then $\Pi(\chi)$ is nonconvex, i.e. the so-called double-well potential, and  the problem (\ref{eq-beamfix}) has three fixed points.
In this case, the beam is in post-buckling state. In order to solve this challenging nonconvex stationary problem, a canonical dual finite element method has been developed recently \cite{ali-gao-amma}, which shows that for a given nontrivial lateral load $q(x)$, the
nonlinear differential equation (\ref{eq-beam}) has three solutions at each coordinate $x\in [0,L]$, one is   a global minimizer of $\Pi(\chi)$, which is corresponding to a globally stable post-buckling state of the beam, one is local minimizer of $\Pi(\chi)$, which is corresponding to a locally stable post-buckling state, and one local minimizer of $\Pi(\chi)$, which is corresponding to the unstable post-buckling state.
The numerical results shown that the locally stable post-buckling configuration is extremely sensitive to the
external load $q(x)$ and numerical precision used in the program.

For unilateral post-buckling problems, the feasible set $\calX_a$ has usually inequality constraints.
For example, a  simply supported beam on rigid foundation subjected to a  downward lateral load $q(x)   \;\;\forall x\in[0,L]$,
this feasible set is a convex cone:
\[
\calX_a = \{ \chi(x) \in C^2[0,L]| \;\; \chi(x) \ge 0 \;\; \forall x\in [0, L], \;\; \chi(0)=\chi(L) = 0 \}.
\]
Due to the inequality constraint in $\calX_a$, the stationary condition
of the problem (\ref{eq-beamfix}) leads to not only a so-called variational inequality \cite{liu-gao-wang}
\eb
\int_0^L  (\ww - \chi)  \delta \Pi(\chi) \dx  \ge 0   \;\; \forall \ww(x) \in \calX_a,
 \ee
 where $\delta\Pi(\chi) =  \chi_{xxxx} -   \frac{3}{2}  \alp \chi_x^2 \chi_{xx}  + \lam \chi_{xx} -  q$ is the \G derivative of $\Pi(\chi)$,
but also the well-known complementarity condition
 \eb
(  \chi_{xxxx} -   \frac{3}{2}   \alp\chi_x^2 \chi_{xx}  + \lam \chi_{xx} -  q) \chi(x) = 0 \;\; \forall x \in [0,L]
\ee
Since the   contact region (i.e. on which $\chi(x) = 0$) remains unknown till the problem is solved,
 the problem (\ref{eq-beamfix}) is the combination of the nonlinear free-boundary value problem, non-monotone  variational inequality, and the nonconvex variational analysis. This problem could be one of  the most challenging problems in  nonconvex analysis, which deserves seriously study in the future.
}
\end{exam}
Canonical duality-triality is a  methodological  theory  which can be used not only for modeling complex systems within a unified framework, but also for solving real-world  problems with a unified methodology.
This theory was developed  originally from Gao and Strang's work in nonconvex mechanics
 \cite{gao-strang1989} and has been applied successfully for solving a large class of challenging problems
 in both nonconvex analysis/mechancis and global optimization, such as
 chaotic dynamics \cite{lat-gao-chaos, ruan-gao-ima},
 phase transitions in solids \cite{gao-yu}, post-buckling of large deformed beam \cite{ali-gao-amma},
nonconvex and discrete optimization
 \cite{gao-chen2014,gao2005,gao2007,gao-ruan2009,gao-ruan2010,wang-fang-gao2008}.
A comprehensive review on this theory and breakthrough
 from recent challenges are given in \cite{gao-lat-ruan-amma}.

The goal of this paper is to apply the canonical duality  theory for solving the challenging fixed point problem.
The rest of this paper is arranged as follows.
Based on the concept of objectivity, a fixed point problem and its  canonical dual
  are formulated in the next section.
  Analytical solutions for  a general fixed point problem with
   sum of exponential functions and nonconvex polynomial are discussed in Sections 3.
    Analytical solutions for  a general fixed point problem with
   sum of Logarithmic and quadratic functions  are discussed in Sections 4.
  Special examples  are illustrated  in Sections 5 and 6.
The paper is ended with conclusions and future work.

\section{Canonical Dual Solutions to Fixed Point Problems}
According to the canonical duality,  the linear measure $\epsilon = D \bx$ can't be used directly  for
studying duality relation due to the objectivity.
Also,  the linear operator can't change the nonconvexity of $W(D\bx)$.
We first introduce the canonical transformation.

\begin{definition}[Canonical Function and Canonical Transformation]$\;$ \hfill

A real-valued function $\VV:\calE_a \rightarrow \real$ is called canonical if the duality mapping
$ \partial \VV: \calE_a \rightarrow \calE_a^*$ is one-to-one and onto.

For a given nonconvex function $W:\calW_a \rightarrow \real$, if there exists a geometrically admissible mapping
$\Lam:\calW_a \rightarrow \calE_a$ and a  canonical  function $\VV:\calE_a \rightarrow \real$ such that
\eb
W(\bepsilon) = \VV ( \Lam(\bepsilon)),  \label{eq-ct}
\ee
then, the transformation (\ref{eq-ct}) is called the canonical transformation % $\Phi(\bxi)$ is called the canonical function,
 and  $\bxi = \Lam(\bepsilon)  $ is called the canonical  measure.
\end{definition}

By this definition, the  one-to-one duality relation
$\bvsig = \partial \VV(\bxi) : \calE_a \rightarrow
\calE^*_a$  implies that the  canonical function $ \VV(\bxi)$  is differentiable
 and its  conjugate function $ \VV^*:\calE^*_a \rightarrow \real$
 can be uniquely defined by the Legendre transformation \cite{gao-book00}
\eb
 \VV^*(\bvsig) = \{ \la \bxi ; \bvsig \ra -  \VV(\bxi) | \; \bvsig = \partial  \VV(\bxi) \} ,
\ee
where $\la \bxi ; \bvsig \ra $ represents the bilinear form on $\calE$ and its dual space $\calE^*$.
In this case,  $\VV:\calE_a \rightarrow \real$ is a canonical function if and only if
  the following canonical duality relations hold on $\calE_a \times \calE^*_a$:
 \eb
  \bvsig = \partial  \VV(\bxi)  \;\; \Leftrightarrow \;\;   \bxi= \partial \VV^*(\bvsig ) \;\; \Leftrightarrow \;\;
 \VV(\bxi) +  \VV^*(\bvsig) = \la \bxi ; \bvsig \ra. \label{eq-cdr}
  \ee

Let $Q(\bx)=\frac{1}{2}\|\bx\|^2 + \bx^T \bff$. Replacing $ \VV(\Lambda(\bx)) $ in  the target function $\Pi(\bx)$  by the Fenchel-Young equality
 $\VV(\bxi)  =  \la \bxi ; \bvsig \ra - \VV^*(\bvsig)$,
  the  total complementary function  (see \cite{gao-jogo00})
   $\Xi: \calX_a \times \calE^*_a \rightarrow \bbR $  can be obtained as
 \beq
 \Xi(\bx, \bvsig ) = \la \Lambda(\bx)  ; \bvsig \ra - \VV^*(\bvsig) - Q(\bx).
 \eeq
 By this total complementary function, the canonical dual of $\Pi(\bx)$ can be obtained as
 \beq
 \Pi^d(\bvsig) =  \inf \{ \Xi(\bx, \bvsig) | \;\; \bx \in \calX \} = Q^{\Lambda}(\bvsig) - \VV^*(\bvsig),
 \eeq
 where $Q^\Lambda:\calE^*_a\rightarrow \bbR \cup\{ - \infty\}$ is the so-called $\Lambda$-conjugate of $Q(\bx)$ defined by
 (see \cite{gao-jogo00})
 \eb
 Q^\Lambda(\bvsig) =  \sta \{ \la \Lambda(\bx) ;  \bvsig \ra - Q(\bx) \; | \;\; \bx \in \calX \}.
 \ee
 Let $\calS_a \subset  \calE^*_a$ be an admissible set such that on which,  $Q^\Lam(\bvsig)$ is well-defined.
 If $\Lam(\bx)$ is a homogeneous quadratic operator, i.e. $\Lam(\alp \bx) = \alp^2 \Lam(\bx)$,  then the total complementary function
 \beq
 \Xi(\bx, \bvsig ) = \half \la  \bx,   \bG(\bvsig )   \bx \ra - \VV^*(\bvsig)  - \la \bx , \bff \ra  ,
 \eeq
 where $\bG(\bvsig) = \bH(\bvsig) - \bI$,   $\bH(\bvsig) = \nabla^2_{\bx}  \la \Lam(\bx); \bvsig \ra$, and $\bI$ is an identity matrix  in $\calX$.
 In this case,  the $\Lam$-conjugate
 $Q^\Lam$ is simply defined by
 \eb
  Q^\Lambda(\bvsig) = -\half  \la   \bG^{-1} (\bvsig)   \bff, \bff \ra  ,
  \ee
 and $\calS_a = \{ \bvsig \in \calE^*_a | \;\; \det \bG(\bvsig) \neq 0 \}$.
 Thus, the canonical dual problem $(\calP^d)$ can be proposed in the following:
 \eb
 (\calP^d): \;\; \barbvsig = \arg \sta \left\{ \Pi^d(\bvsig) | \;\; \bvsig \in \calS_a \right\}. \label{eq-cdfix}
 \ee
 By the canonical duality theory,  we have the following results.
\begin{thm}[Analytic Solution and Complementary-Dual Principle]\label{th:AnalSolu} $\;$\newline
For  a given $\bff$, if $\barbvsig \in \calS_a$ is a solution to $(\calP^d)$, then
\beq\label{eq:solvedx}
\barbx= \bG^{-1}(\barbvsig)\bff
\eeq
is a  solution to the problem $(\calP)$  and
\beq\label{eq:nogap}
\Pi(\barbx)=\Xi(\barbx,\barbvsig)=\Pi^d(\barbvsig).
\eeq
If $F(\bx)$ is a potential operator, then $\barbx$ is also a solution to the  fixed point problem $(\calP_0)$.
If $F(\bx)$ is a non-potential operator, then $\barbx$ is  a solution to the  fixed point problem $(\calP_0)$ only if $\barbx$ is
a global minimizer of $\Pi(\bx)$.
\end{thm}
{\bf Proof}. By the canonical duality theory we know that $(\barbx, \barbvsig)$ is a critical point of $\Xi(\bx, \bvsig)$ if and only if $\barbx$ is a critical point of $\Pi(\bx)$ and $\barbvsig$ is a critical point of $\Pi^d(\bvsig)$.
It is easy to prove that the criticality condition $\nabla \Xi(\barbx, \barbvsig) = 0$ leads to the following canonical equations:
\eb
\bG(\barbvsig) \barbx = \bff ,\;\;\;
\Lam(\barbx) = \partial \VV^*(\barbvsig).
\ee
The first equation is the canonical equilibrium equation, which leads to the analytical solution (\ref{eq:solvedx}).
By the canonical duality, the canonical duality equation $\Lam(\barbx) = \partial \VV^*(\barbvsig)$ leads to the complementary-duality relation
(\ref{eq:nogap}).  The theorem is proved by Lemma 1. \hfill $\Box$ \\

Theorem \ref{th:AnalSolu} shows that the solution to the fixed point problem  depends analytically on
the canonical dual solution, and there is no duality gap between the primal
problem ($\calP$) and the canonical dual problem ($\calP^d$).
By the fact  that the problem $(\calP)$ may have many fixed points,
in order to identify the extremality of these fixed points,
 we introduce two open sets:
\begin{eqnarray*}
\calS_a^+ &=&\{\bvsig \in \calS_a | \bG(\bvsig)\succ 0\},\\
\calS_a^- &=&\{\bvsig \in \calS_a| \bG(\bvsig)\prec 0\}
\end{eqnarray*}
where $\bG(\bvsig)\succ 0$ means that $\bG$ is a positive definite matrix,
and $\bG(\bvsig) \prec 0$ means that $\bG$ is a negative definite matrix.
Also according to the terminology used in the canonical duality theory, a  neighborhood of a critical point is an open set
containing only one critical point.

\begin{thm}[Triality Theorem]\label{th:triality}
 Suppose that $\barbvsig$ is a solution to $(\calP^d)$
  and $\barbx=\bG^{-1}(\barbvsig)\bff$.
If $\barbvsig \in \calS_a^+$,  then  $\barbx$ is a globally stable fixed point and
\begin{eqnarray}\label{statement1}
\Pi(\bar{\bx})   =\min_{\bx\in \calX_a } \Pi(\bx)=\max_{\bvsig\in\calS_a^+} \Pi^d(\bvsig)=\Pi^d(\bar{\bvsig}).
\end{eqnarray}

\item If $\barbvsig\in\calS_a^-$, then
 $\barbx$ is a local maximizer of $\Pi(\bx)$ iff $\barbvsig\in\calS_a^-$
  is a local maximizer of $\Pi^d$ and   on the  neighborhood $\calX_o \times \calS_o \subset \calX_a \times \calS_a^- $  of $(\barbx,\barbvsig)$,
   we have
    \begin{equation}\label{statement2}
    \Pi(\barbx)   =\max_{\bx\in\calX_o} \Pi(\bx)=\max_{\bvsig\in\calS_o} \Pi^d(\bvsig)=\Pi^d(\barbvsig).
    \end{equation}
 Moreover,  $\barbx$ is a locally unstable fixed point of  $F$   if it  is a potential operator.

\item If $\barbvsig\in\calS_a^-$  and $\dim \calX_a = \dim \calS_a$, then  $\barbx$ is a  local minimizer iff $\barbvsig \in \calS^-_a$ is a local minimizer of $\Pi^d(\bvsig)$ and on the  neighborhood $\calX_o \times \calS_o \subset \calX_a \times \calS_a^- $ ,
    \begin{equation}\label{statement3}
    \Pi(\barbx)   =\min_{\bx\in\calX_o} \Pi(\bx)=\min_{\bvsig\in\calS_o} \Pi^d(\bvsig)=\Pi^d(\barbvsig) .
    \end{equation}
 Moreover,  $\barbx$ is a locally stable fixed point of   $F$   if it  is a potential operator.
 \end{thm}

%Details proof please refer to \cite{gao-ruan-latorre2015}.
This theorem is an application of the triality theory. Detailed proof can be found in \cite{gao-wu-amma}.
The statement (\ref{statement1}) is   the so-called {\em canonical min-max duality}. This statement shows that
the global stable fixed point problem n  is equivalent to
a concave maximization problem
\eb
(\calP^\sharp): \;\; \max \{ \Pi^d(\bvsig) | \;\; \bvsig \in \calS^+_a \}.
\ee
Since  the feasible space $\calS^+_a$ is  an open convex set, this canonical dual problem     can be solved easily by well-developed nonlinear  optimization techniques.
The second statement (\ref{statement2}) is the {\em canonical double-max duality} and
the third statement  (\ref{statement3}) is the  {\em canonical double-min duality}.
  For potential operator $F$, these two statements can be used to identify locally unstable and stable  fixed points, respectively.

\section{Application to Exponential and Polynomial Functions}
As our first application,   the objective function is assumed to be
\begin{eqnarray*}
\WW(\DD \bx)=\alpha \exp(\frac{1}{2}\|\DD_1 \bx \|^2)
+\frac{1}{2} \beta (\frac{1}{2} \|\DD_2 \bx\|^2 -\lambda)^2,
\end{eqnarray*}
where $\DD_1  \in \mathbb{R}^{m\times n}$ and $ \DD_2 \in \mathbb{R}^{p\times n}$ are two given matrices,
$\alpha$, $\beta$, $\lambda$ are real numbers.
Clearly,  for a given $\lam > 0$,  $\WW(\DD \bx)$ is nonconvex and
\begin{eqnarray*}
F(\bx)=\nabla P(\bx)=\alpha \exp(\frac{1}{2}\|\DD_1 \bx \|^2)(\DD_1^T\DD_1)\bx
+\beta (\frac{1}{2} \|\DD_2 \bx\|^2 -\lambda)(\DD_2^T\DD_2)\bx-\bff
\end{eqnarray*}
is a non-monotone operator.
In this case, the fixed point problem $(\calP_0)$ can be equivalently written as
\begin{eqnarray*}
\barbx = \arg \sta \left\{ \Pi(\bx)=\alpha \exp(\frac{1}{2}\|\DD_1 \bx \|^2)
+\frac{1}{2} \beta (\frac{1}{2} \|\DD_2 \bx\|^2 -\lambda)^2-\frac{1}{2}\|\bx\|^2-\bx^T\bff \;\; | \;\; \bx \in \real^n \right\} .
\end{eqnarray*}
Clearly, traditional methods for solving this   nonlinear fixed point problem in $\real^n$  are difficult.
However, by the canonical duality theory, this problem can be solved easily in $\real^2$.

The canonical measure in this problem can be given as
\begin{eqnarray*}
\bxi = \Lam (\bx)=
\left(\begin{array}{c}
\xi_{1}  \\
\xi_{2}
\end{array}
\right)
=\left(\begin{array}{c}
 \frac{1}{2}\|\DD_1\bx\|^2\\
 \frac{1}{2}\|\DD_2\bx\|^2
\end{array}
\right)  \; :  \real^n \rightarrow   \calE_a = \{ \bxi   \in \real^2| \;\; \xi_1, \xi_2 \ge 0 \}.
\end{eqnarray*}
Correspondingly, the canonical function is
\begin{eqnarray*}
V(\bxi)
=
\left(\begin{array}{c}
V_1(\xi_1)  \\
V_2(\xi_2)
\end{array}
\right)
 = \left(\begin{array}{c}
 \alpha \exp(\xi_1) \\
  \frac{1}{2} \beta (\xi_2-\lambda)^2
  \end{array} \right)  ,
\end{eqnarray*}
and the canonical dual variable is
\begin{eqnarray*}
\bvsig =
\left(\begin{array}{c}
\vsig_{1}  \\
\vsig_{2}
\end{array}
\right)
=\left(\begin{array}{c}
 \nabla V_1(\xi_{1})\\
 \nabla V_2(\xi_{2})
\end{array}
\right)
=\left(\begin{array}{c}
\alpha \exp(\xi_1)\\
 \beta(\xi_2-\lambda)
\end{array}
\right) \;\; :\calE_a \rightarrow \calE^*_a = \{ \bvsig \in \real^2 | \;\; \vsig_1 \ge 0 \}.
\end{eqnarray*}
By the Legendre transformation, the conjugate function
  $V^*(\bvsig)$  is uniquely  defined as
\begin{eqnarray*}
V^*(\bvsig)
=\left(\begin{array}{c}
V_1^*(\vsig_1)\\
V_2^*(\vsig_2)
\end{array}
\right)
  =    \left(\begin{array}{c}
  (\ln(\vsig_{1}/\alpha)-1) \vsig_1 \\
  \frac{1}{2 \beta} \vsig_2^2 + \lambda \vsig_2 \end{array}  \right).
\end{eqnarray*}
 Since the canonical measure in this application is a homogeneous quadratic operator,
the total complementary function
$\Xi : \mathbb{R}^n\times \mcal{E}_a^*\rightarrow \mathbb{R}$ has the following form
\begin{eqnarray*}
\Xi(\bx,\bvsig)=\frac{1}{2}\bx^T \bG(\bvsig)\bx-\bx^T \bff
-(\ln(\vsig_1/\alpha)-1)\vsig_1-(\frac{1}{2 \beta}\vsig_2^2 + \lambda \vsig_2),
\end{eqnarray*}
where,
\begin{eqnarray*}
\bG(\bvsig)=\vsig_1 \DD_1^T \DD_1 +\vsig_2 \DD_2^T \DD_2 - \bI.
\end{eqnarray*}
On the canonical dual feasible space $
\calS_a= \{ \bvsig=[\vsig_1,\vsig_2]^T \in \real^2 | \det( \bG(\bvsig)) \neq 0 \},$
  the canonical
dual problem can be formulated as
\eb
(\calP^d):\;\; \barbvsig = \arg \sta \left\{
P^d(\bvsig)=
-\half \bff^T \bG^{-1}(\bvsig) \bff  -(\ln(\vsig_1/\alpha)-1)\vsig_1-(\frac{1}{2 \beta}\vsig_2^2 + \lambda \vsig_2) \; |
\;\;\bvsig \in \calS_a \right\}.
\ee

\subsection*{Example 1}
\begin{figure}[h!]
\centering
\mbox{\resizebox{!}{2.0in}{\includegraphics{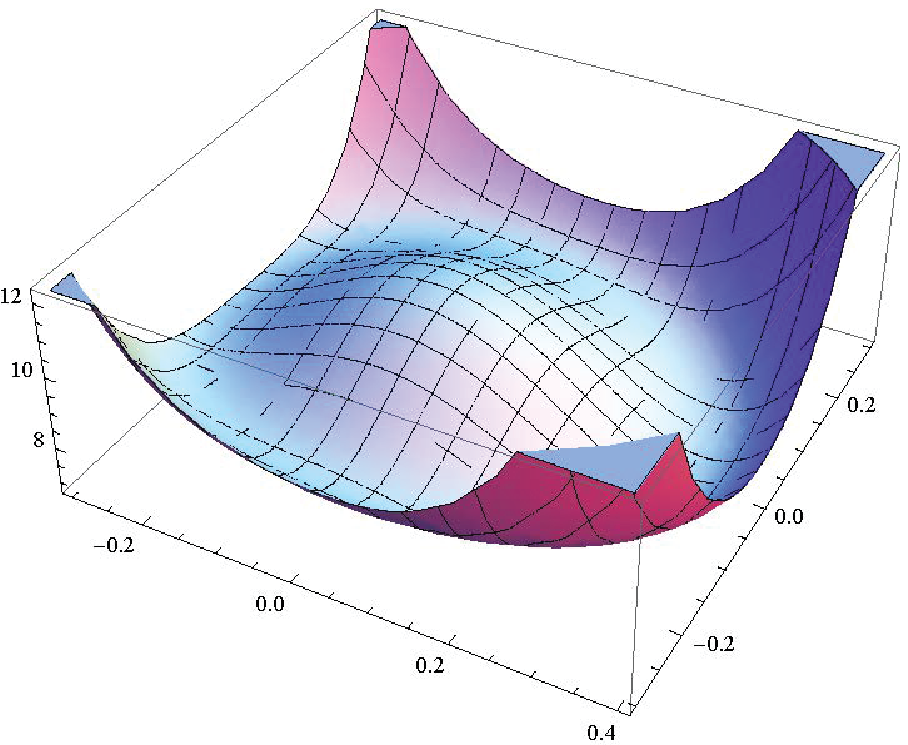}}\quad\quad
\resizebox{!}{2.0in}{\includegraphics{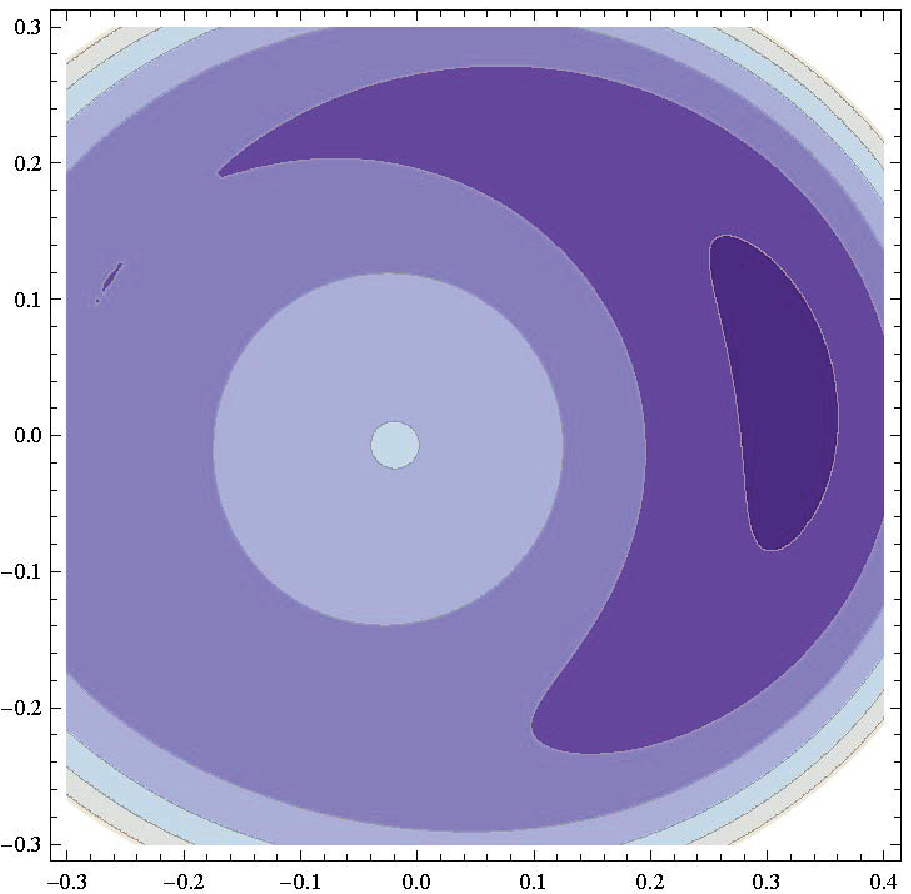}}}\vspace{-0.5cm}
\caption{\label{figure-ex1} Graphs of $\Pi(x_1,x_2)$ and its contour for Example 1.}\label{fig-a2}
\end{figure}

Let $n=2$, $\alpha=6$, $\beta=8$, $\lambda=1$ and
\begin{eqnarray*}
\DD_1=\left[
\begin{array}{cc}
2&0\\
0&3
\end{array}
\right],
\DD_1=\left[
\begin{array}{cc}
4&0\\
0&5
\end{array}
\right], \;\;
\bff  =\left[
\begin{array}{c}
2\\
1
\end{array}
\right],
\end{eqnarray*}
then the primal function (see Figure \ref{figure-ex1})
\begin{eqnarray*}
 \Pi(x_1,x_2)=
6 \exp(2x_1^2+4.5x_2^2)+4(8x_1^2+12.5x_2^2-1)^2-\frac{1}{2}(x_1^2+x_2^2)-2x_1-x_2.
\end{eqnarray*}
The corresponding canonical dual function is
\begin{eqnarray*}
\Pi^d(\vsig_1,\vsig_2)=-\frac{1}{2}
(\frac{4}{4\vsig_1+16\vsig_2-1}+\frac{1}{9\vsig_1+25\vsig_2-1})
-\vsig_1 (\ln(\vsig_1/6)- 1) -(\frac{1}{16}\vsig_2^2+\vsig_2)
\end{eqnarray*}
\begin{figure}[h!]
\centering
\mbox{\resizebox{!}{2.0in}{\includegraphics{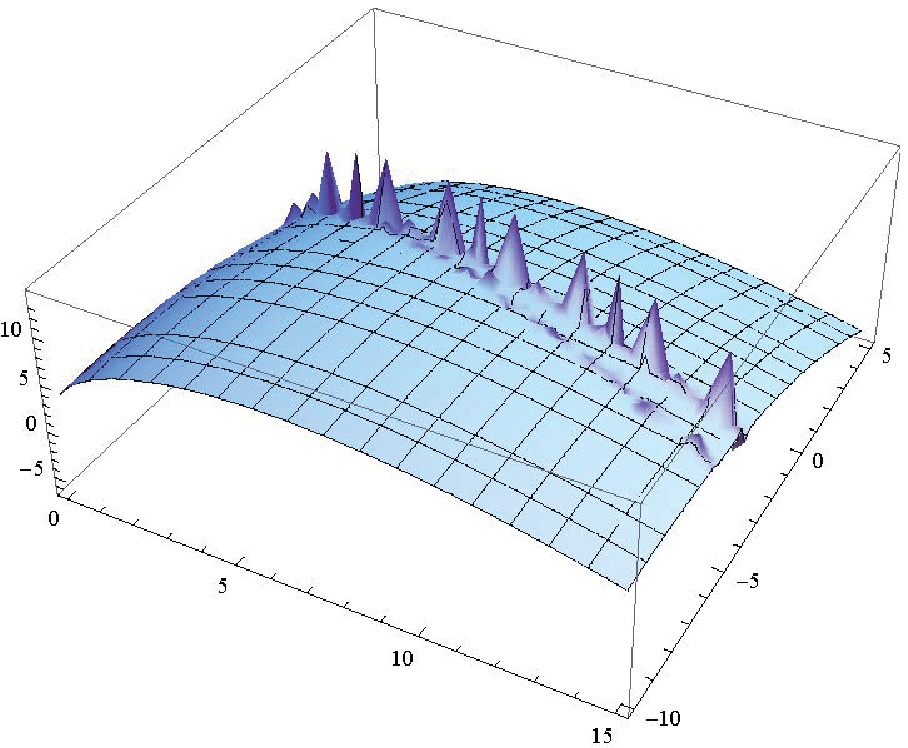}}\quad\quad
\resizebox{!}{2.0in}{\includegraphics{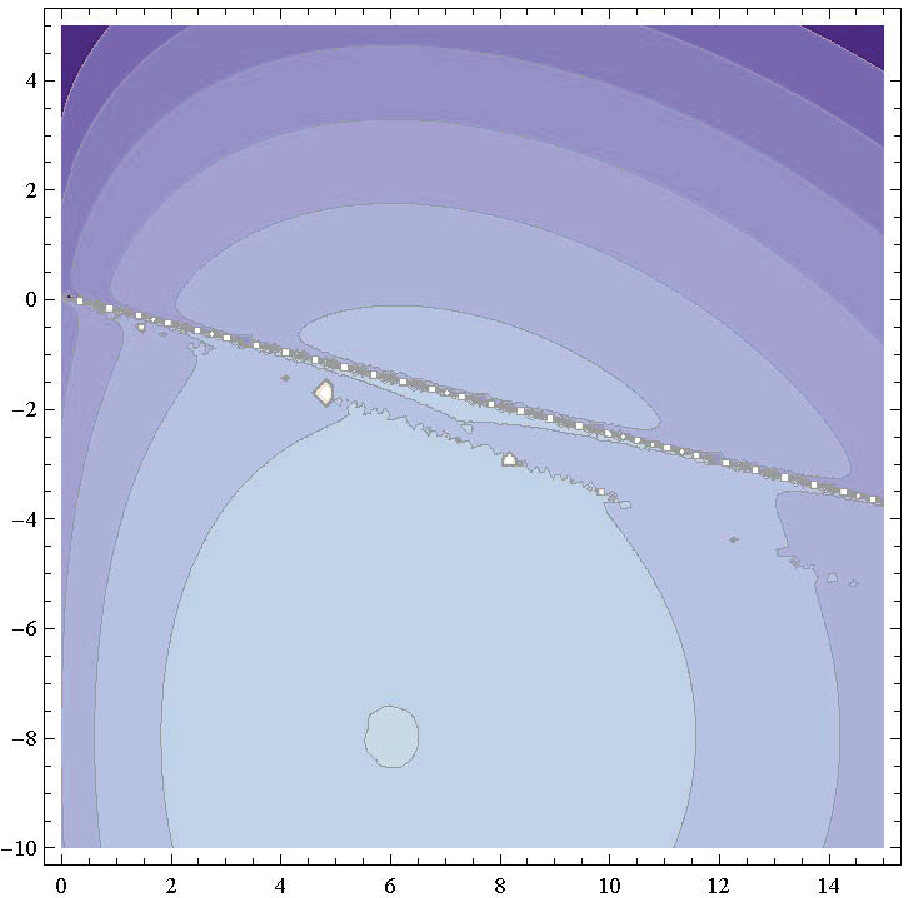}}}\vspace{-0.5cm}
\caption{\label{figure-dual-ex1} Graphs of $\Pi^d(\vsig_1, \vsig_2)$ and its contour for Example 1.}
\end{figure}
Its graph is shown by  Figure \ref{figure-dual-ex1}.
It is easy to find that the canonical dual problem $(\calP^d) $ has three solutions:
\begin{eqnarray*}
& & \bvsig^1  =[7.38697,-1.39206]^T \in \calS_a^+, \\
& & \bvsig^2  =[6.00566,-7.97189]^T\in \calS_a^- , \\
& & \bvsig^3  =[7.3106,-2.23695]^T\in \calS_a^-.
\end{eqnarray*}
By Theorem \ref{th:AnalSolu}
we have three primal solutions:
\begin{eqnarray*}
& & \bx^1 =[0.318731, 0.0325932]^T, \\
& & \bx^2 =[-0.0191337,-0.00683777]^T, \\
& & \bx^3 =[-0.264945, 0.112718]^T .
\end{eqnarray*}
It is easily to check that
\begin{eqnarray*}
& & \Pi(\bx^1)=\Pi^d(\bvsig^1)=6.78671 , \\
& & \Pi(\bx^2)=\Pi^d(\bvsig^2)=10.0225 ,\\
& & \Pi(\bx^3)=\Pi^d(\bvsig^3)=7.99906.
\end{eqnarray*}
By Theorem \ref{th:triality}
we know that $\bx^1 $ is a global minimizer of $\Pi(\bx)$,
 $\bx^2$ is  a local maximizer  of $\Pi(\bx)$,   and
$\bx^3=[-0.264945, 0.112718]^T$ is a local minimizer of $\Pi(\bx)$ (see Figure \ref{figure-ex1}).
By the fact that
\begin{eqnarray*}
x_1^i &=& F_1(x_1^i,x_2^i)=6\exp(2x_1^i+4.5 x_2^i)4x_1^i+8(8x_1^i+12.5x_2^i-1)16x_1^i-2,\\
x_2^i &=& F_2(x_1^i,x_2^i)=6\exp(2x_1^i+4.5 x_2^i)9x_2^i+8(8x_1^i+12.5x_2^i-1)25x_2^i-1,
\end{eqnarray*}
hold for all $i=1,2,3$, we know that $\{ \bx^i\} $ ($i=1,2,3$)  are all fixed points.

\section{ Application to Logarithmic and Quadratic Function}
In this application, we let
\begin{eqnarray*}
\WW(\DD \bx)=c_1 \|\DD \bx\|^2+c_2\|D\bx\|^2 \log \|\DD \bx\|^2
\end{eqnarray*}
where $\DD \in \mathbb{R}^{m\times n}$ is matrix,
$c_1$, $c_2$ are real numbers.
Clearly, $\WW(\DD \bx)$ is nonconvex and
\begin{eqnarray*}
F(\bx)=\nabla P(\bx)
=2c_1(\DD^T \DD)\bx+2c_2((\DD^T \DD) \bx) (\log\|\DD \bx\|^2+1)
\end{eqnarray*}
is non-monotone. The fixed point problem $\bx = F(\bx)$ can be reformulated as
\begin{eqnarray*}
(\calP): \;\; \barbx = \arg \sta \left \{ \Pi(\bx)=c_1 \|\DD \bx \|^2+c_2\|\DD \bx \|^2 \log \|\DD \bx\|^2
-\frac{1}{2}\|\bx\|^2-\bx^T\bff | \;\; \bx \in \real^n \right\}.
\end{eqnarray*}

By using the canonical measure
\[
\xi = \Lam (\bx)= \|\DD\bx\|^2 : \; \real^n \rightarrow  \calE_a = \real^+= \{ \xi \in \real| \; \xi \ge 0 \},
\]
the  canonical function is  $V(\xi)=c_1 \xi+c_2 \xi (\log \xi +1)$ and its Legendre conjugate is
 \begin{eqnarray*}
V^*(\vsig)=c_2 \exp [\frac{1}{c_2}(\vsig-c_1)-1],
\end{eqnarray*}
which is convex on its domain $\calE^*_a = \real $.
 In this case, we have  the total complementary function
\begin{eqnarray*}
\Xi(\bx,\vsig)=\frac{1}{2}\bx^T \bG(\vsig)\bx-\bx^T \bff
-c_2 \exp [\frac{1}{c_2}(\vsig-c_1)-1],
\end{eqnarray*}
where $
\bG(\vsig)=2\vsig\DD^T \DD  - \bI$ and the canonical dual problem is
 \eb
 (\calP^d): \;\; \barvsig = \arg\sta \left \{
P^d(\vsig)=
-\half \bff^T \bG^{-1}(\vsig) \bff  -c_2 \exp [\frac{1}{c_2}(\vsig-c_1)-1]\;\; | \;\; \vsig \neq 0 \right\} .
\ee
 \subsection*{Example 2}
\begin{figure}[h!]
\centering
\mbox{\resizebox{!}{2.0in}{\includegraphics{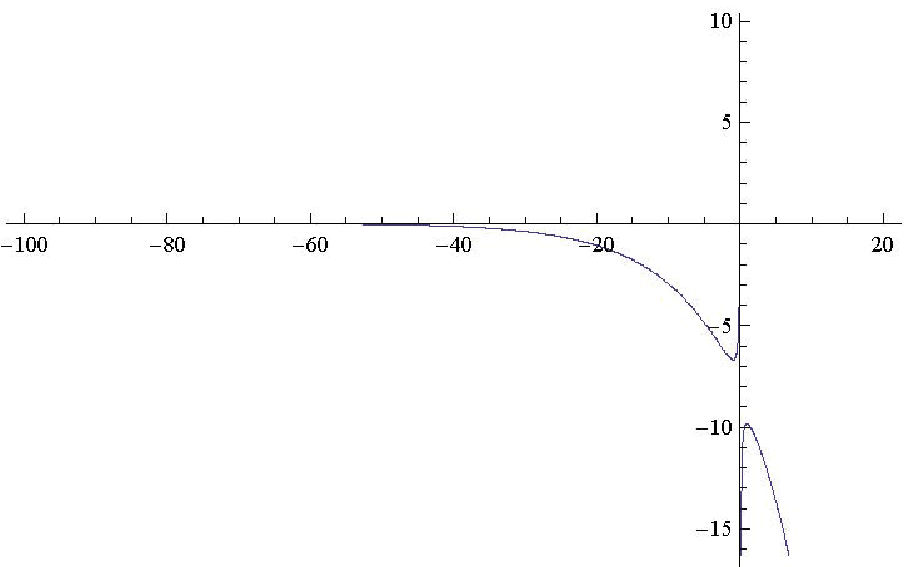}}}
\vspace{-0.5cm}
\caption{\label{figure-dual-ex2} Graph  of $\Pi^d(\vsig)$ for Example 2.}
\end{figure}

We first let $m=n=2$, $c_1=-8$, $c_2=10$,  and
\begin{eqnarray*}
\DD=\left[
\begin{array}{cc}
3&0\\
0&4
\end{array}
\right],
\;\; \bff =\left[
\begin{array}{c}
-5\\
2
\end{array}
\right].
\end{eqnarray*}
The  primal function
\begin{figure}[h!]
\centering
\mbox{\resizebox{!}{2.0in}{\includegraphics{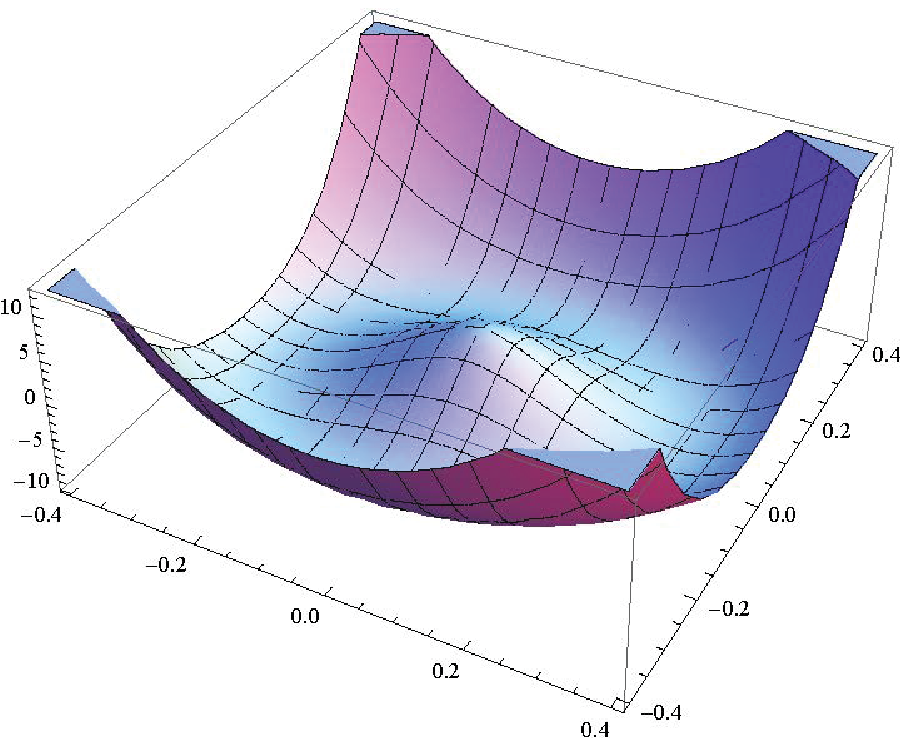}}\quad\quad
\resizebox{!}{2.0in}{\includegraphics{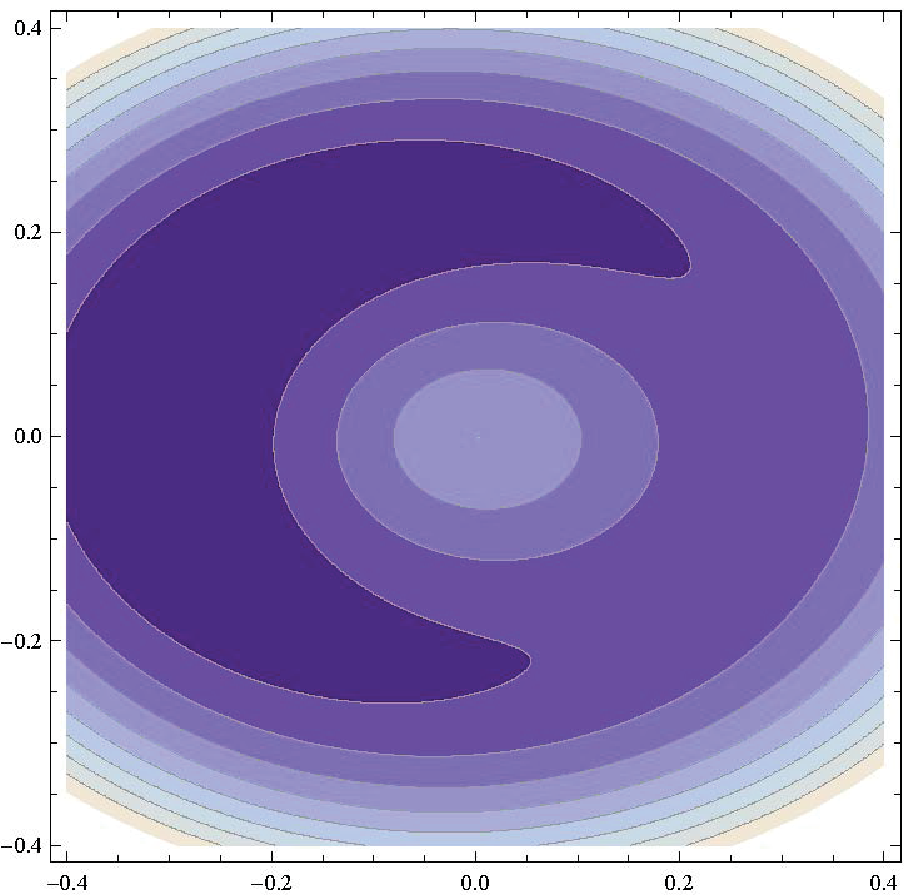}}}\vspace{-0.5cm}
\caption{\label{figure-ex2} Graph  of $\Pi(x_1,x_2)$ and its contour for Example 2.}
\end{figure}
\begin{eqnarray*}
\Pi(x_1,x_2)=
-8(9x_1^2+16x_2^2)+10(9x_1^2+16x_2^2) \log(9x_1^2+16x_2^2)
-\frac{1}{2}(x_1^2+x_2^2)-5x_1+2x_2
\end{eqnarray*}
is nonconvex and its graph is shown in  Figure \ref{figure-ex2}.
The corresponding canonical dual function is
\begin{eqnarray*}
\Pi^d(\vsig)=-\frac{1}{2}
(\frac{25}{18\vsig-1}+\frac{4}{32\vsig-1})
-10\exp[0.1(\vsig+8)-1]
\end{eqnarray*}
For this example, the one dimensional canonical dual problem $(\calP^d)$ can be solved easily (by using Mathematica) to obtain total
three solutions  (see Figure \ref{figure-dual-ex2}):
\begin{eqnarray*}
\vsig^1=0.969642 > \vsig^2  =-0.955077 > \vsig^3 = -91.0174.
\end{eqnarray*}
Correspondingly, the three primal solutions are
\[
\bx^1 = \left [ \begin{array}{c}
-0.303886\\
0.0666033 \end{array}  \right]  , \;\;  \bx^2 = \left [ \begin{array}{c} 0.274855 \\
-0.0633664 \end{array}  \right] , \;\;
\bx^3 =  \left [ \begin{array}{c} 0.00305006 \\
-0.000686446 \end{array} \right] .
\]
It is easy to check that $\bx^i = F(\bx^i), \;\; i=1,2,3$. Therefore, $\{\bx^i\}$ are fixed points.
Since $\vsig^1 \in \calS^+_a = \{ \vsig\in \real| \; \vsig > 0 \} $, we know that $\bx^1$ is a globally stable fixed point.
It is easy to check that $\bx^2 $ is a locally stable fixed point, $\bx^3$ is a locally unstable fixed point, and
\begin{eqnarray*}
\Pi(\bx^1)=\Pi^d(\vsig^1)=-9.84726, \\
\Pi(\bx^2)=\Pi^d(\vsig^2)=-6.69103,  \\
\Pi(\bx^3)=\Pi^d(\vsig^3)=0.00739894 .
\end{eqnarray*}

\subsection*{Example 3}
\begin{figure}[h!]
\centering
\mbox{\resizebox{!}{2.0in}{\includegraphics{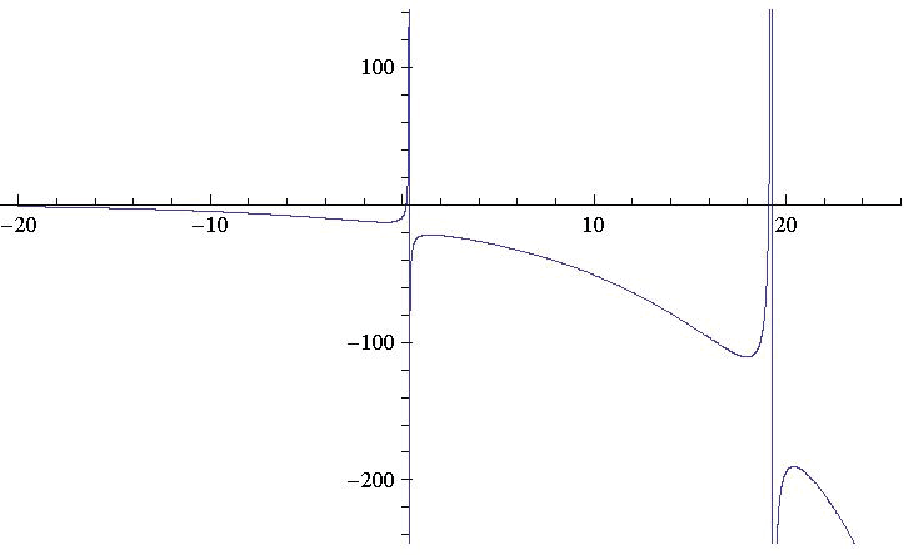}}}
\vspace{-0.5cm}
\caption{\label{figure-dual1-ex3} Graph of $\Pi^d(\vsig)$ for Example 3.}
\end{figure}

\begin{figure}[h!]
\centering
\mbox{\resizebox{!}{2.0in}{\includegraphics{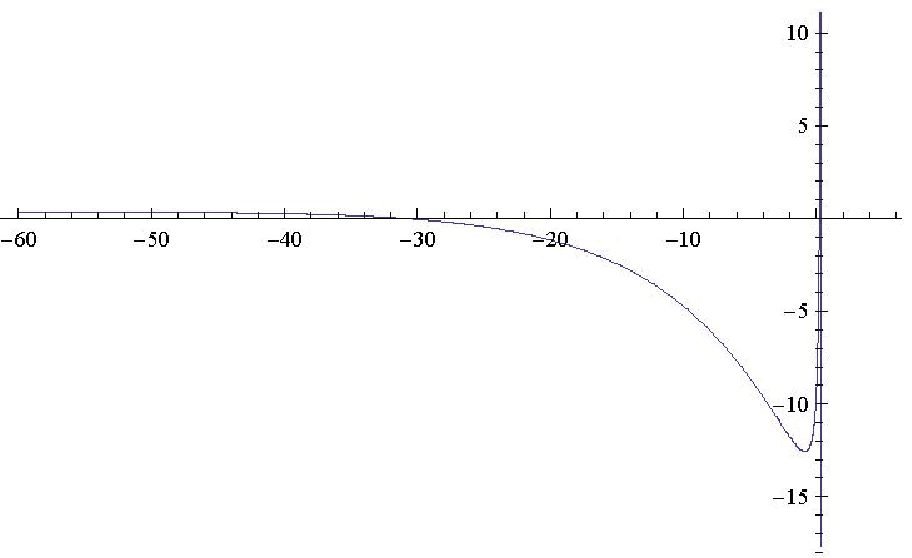}}}
\vspace{-0.5cm}
\caption{\label{figure-dual2-ex3} Graph of $\Pi^d(\vsig)$ for Example 3 around  $\vsig^2$.}
\end{figure}

We now let $m=3$, $n=2$, $c_1=-15$, $c_2=9$,  and
\begin{eqnarray*}
\DD=\left[
\begin{array}{cc}
0.3&0.2\\
0.5&0.6\\
0.4&0.7
\end{array}
\right], \;\; \bff =\left[
\begin{array}{c}
1\\
4
\end{array}
\right],
\end{eqnarray*}
then the primal function is
\begin{figure}[h!]
\centering
\mbox{\resizebox{!}{2.0in}{\includegraphics{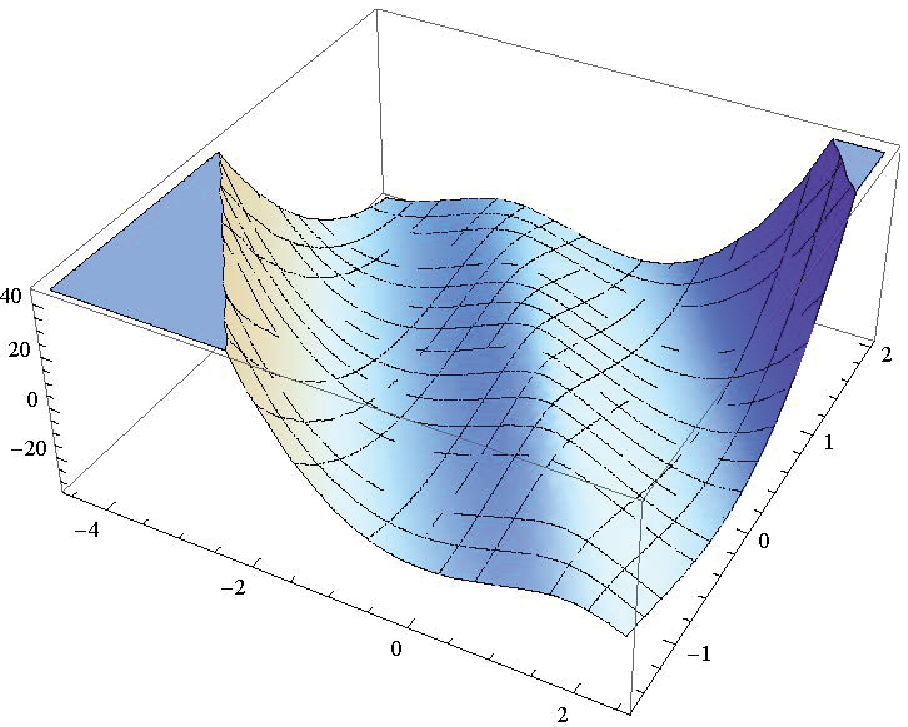}}\quad\quad
\resizebox{!}{2.0in}{\includegraphics{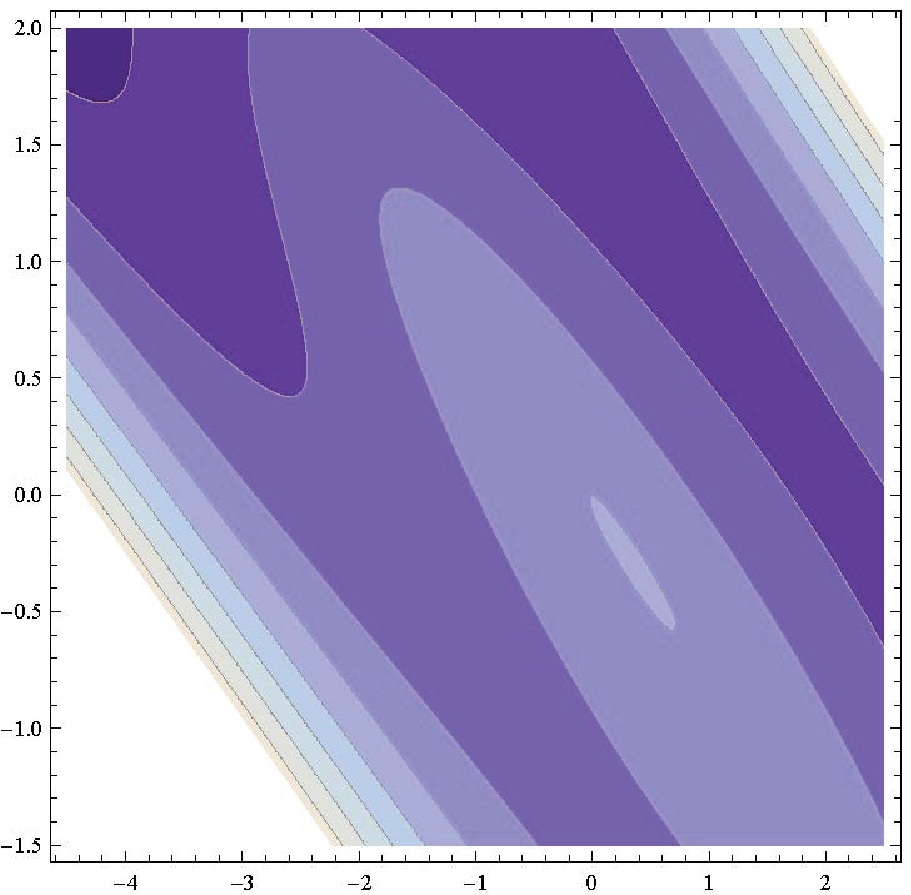}}}\vspace{-0.5cm}
\caption{\label{figure-ex31} Graph of $\Pi(\bx)$ and its contour for Example 3.}
\end{figure}
\begin{figure}[h!]
\centering
\mbox{\resizebox{!}{2.0in}{\includegraphics{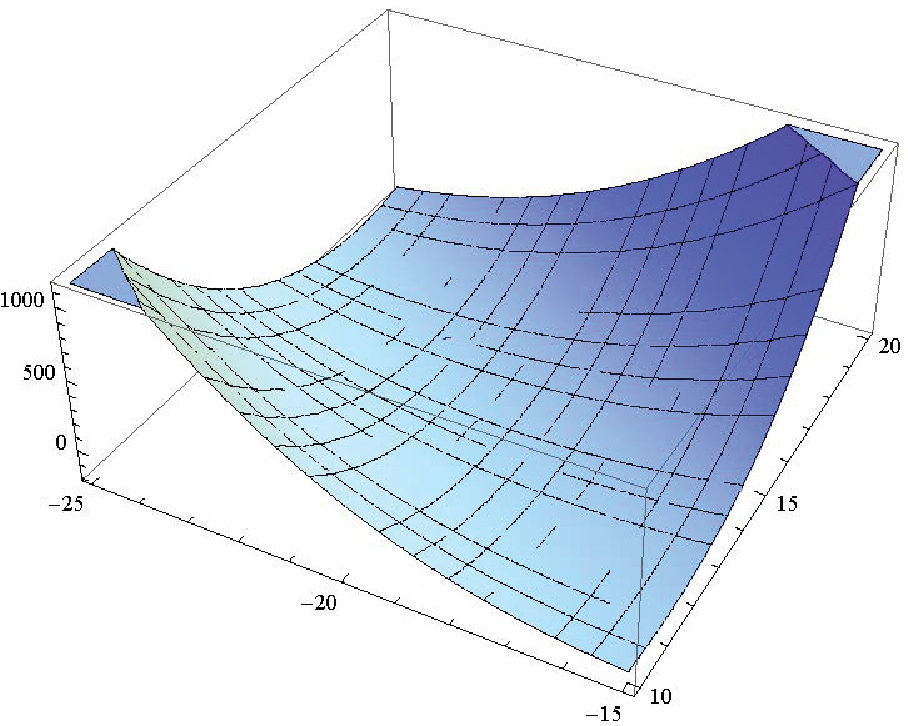}}\quad\quad
\resizebox{!}{2.0in}{\includegraphics{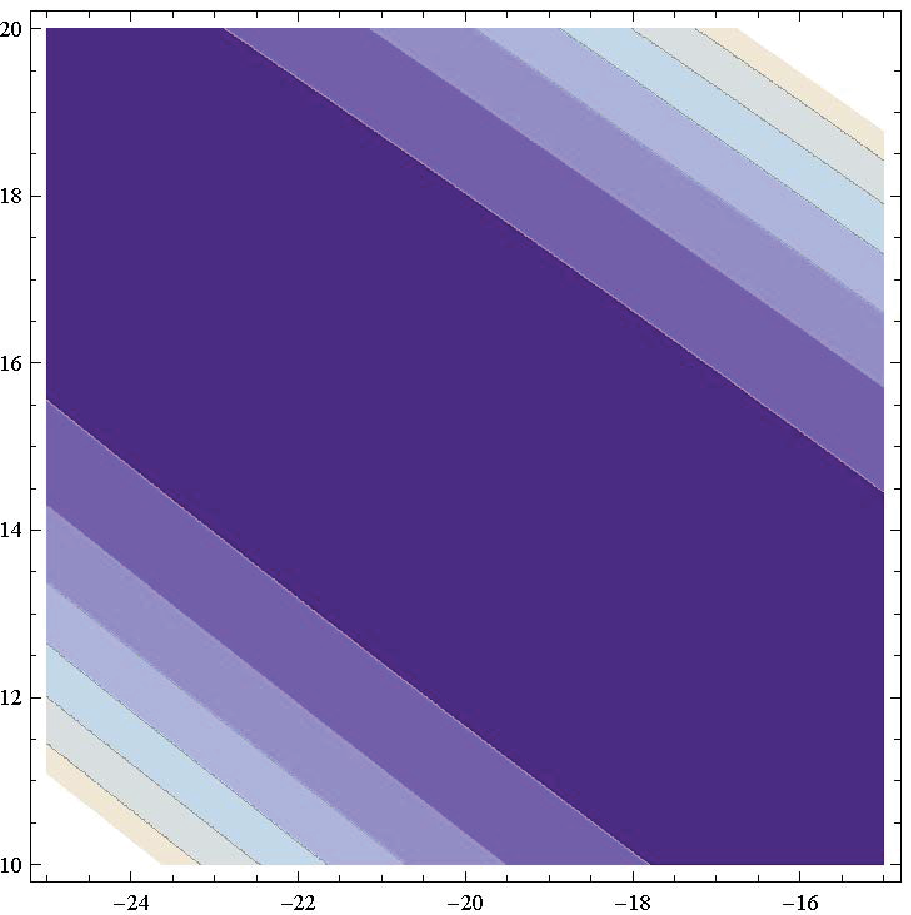}}}\vspace{-0.5cm}
\caption{\label{figure-ex32} Graph  of $\Pi(\bx)$ and its contour for Example 3  around $\bx^1$.}
\end{figure}
\begin{figure}[h!]
\centering
\mbox{\resizebox{!}{2.0in}{\includegraphics{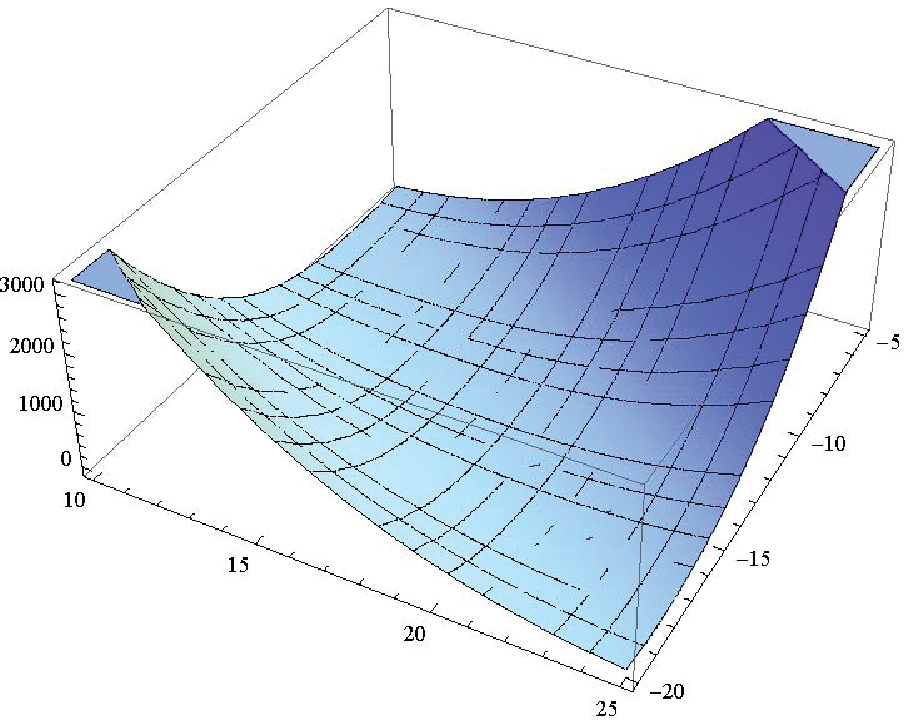}}\quad\quad
\resizebox{!}{2.0in}{\includegraphics{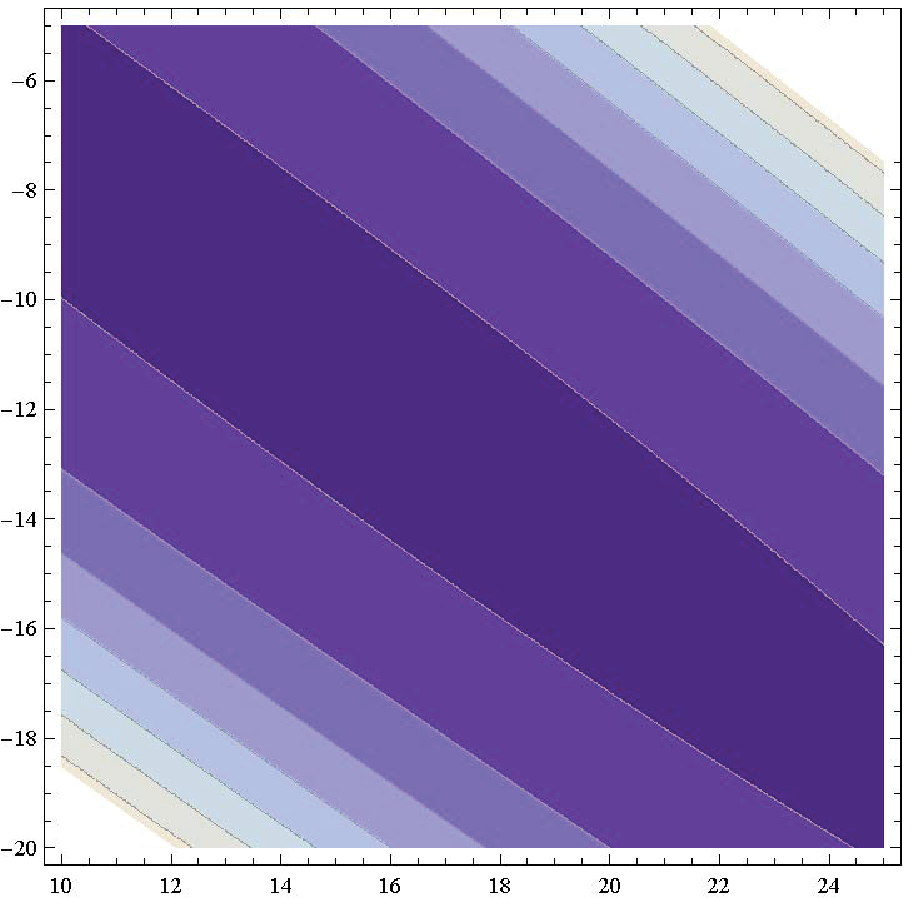}}}\vspace{-0.5cm}
\caption{\label{figure-ex33} Graph of $\Pi(\bx)$ and its contour for Example 3 around  $\bx^2$.}
\end{figure}
\begin{eqnarray*}
\Pi(x_1,x_2) &=&
9(0.5x_1^2+1.28x_1x_2+0.89x_2^2) \log(0.5x_1^2+1.28x_1x_2+0.89x_2^2)\\
&&-15(0.5x_1^2+1.28x_1x_2+0.89x_2^2)-\frac{1}{2}(x_1^2+x_2^2)-x_1-4x_2.
\end{eqnarray*}
Its  graph is a nonconvex surface in $\real^3$, which has multiple critical points, but  their locations can't be find precisely as the surface is rather flat
around these critical points  (see Figure \ref{figure-ex31}-\ref{figure-ex33}).
However, its  canonical dual is a single  valued function
\begin{eqnarray*}
\Pi^d(\vsig)=-\frac{1}{2}
\left[
\begin{array}{cc}
1&4
\end{array}
\right]
\left[
\begin{array}{cc}
\vsig-1&1.28 \vsig\\
1.28 \vsig & 1.78 \vsig -1
\end{array}
\right]^{-1}
\left[
\begin{array}{c}
1\\
4
\end{array}
\right]
-9\exp[\frac{1}{9}(\vsig+15)-1]
\end{eqnarray*}
and from its graph, we can see clearly that it has total five critical points  (see Figure \ref{figure-dual1-ex3}-\ref{figure-dual2-ex3}).
These critical points can be easily obtained by Mathematica:
\begin{eqnarray*}
\vsig^1=20.396 > \vsig^2=17.9735 > \vsig^3  =1.46219 >  \vsig^4=-0.881733 >
\vsig^5=-52.7144 .
\end{eqnarray*}
By Theorem \ref{th:AnalSolu},
we have all the primal solutions:
\[
\bx^1 = \left[ \begin{array}{c}
 -21.57\\
 16.065 \end{array} \right], \;
\bx^2 =  \left[ \begin{array}{c}
18.937 \\
 -13.93 \end{array} \right] , \;
\bx^3 =  \left[ \begin{array}{c}
2.130\\
 0.008 \end{array} \right] ,  \;
\bx^4 =  \left[ \begin{array}{c}
0.546 \\
-1.797 \end{array} \right], \;
\bx^5 =  \left[ \begin{array}{c}
0.323 \\
 -0.272 \end{array} \right].
\]
Since $F(\bx)$ is a potential operator, these stationary points are all fixed points of $F(\bx)$.
It is easy to find that the matrix  $ \bG(\vsig) $ has two singularity points:  $\vsig_1 = 19.266$ and $\vsig_2 = 0.367$, therefore,
\[
  \calS^+_a = \{ \vsig \in \real | \;\; \vsig > 19.266 \}, \;\;
  \calS^-_a = \{ \vsig \in \real | \;\; \vsig < 0.367  \}.
\]
By the facts that $\vsig^1 \in \calS^+_a$ and $\vsig^5\in \calS^-_a$,  we know that  $\bx^1$ is a globally stable fixed point, $\bx^5$ is a locally unstable fixed point.
Although $\vsig^4 \in \calS^-_a$ is a local minimizer of $\Pi^d(\vsig)$, we can't say if $\bx^4$ is a locally stable fixed point since
$\dim \calX_a = 2 \neq \dim \calS_a = 1.$
But by the complementary-dual principle and the order of the canonical dual solutions $\{\vsig^i \}$ , we have
\begin{eqnarray*}
 & & \Pi(\bx^1)=\Pi^d(\vsig^1)=-190.381   \\
 &< &  \Pi(\bx^2)=\Pi^d(\vsig^2)=-110.759 \\
&<&  \Pi(\bx^3)=\Pi^d(\vsig^3)=-21.7036  \\
& < & \Pi(\bx^4)=\Pi^d(\vsig^4)=-12.5735 \\
&<& \Pi(\bx^5)=\Pi^d(\vsig^5)=0.332915.
\end{eqnarray*}

\section{Conclusions}
Based on the canonical duality theory, a unified model is proposed such that the general fixed point problems can be reformulated as a global optimization problem. This model  is directly related to many other challenging problems in variational inequality,
d.c. programming, chaotic dynamics, nonconvex analysis/PDEs, post-buckling of large deformed structures, phase transitions in solids,
and computer science, etc (see \cite{gao-lat-ruan-amma} and references cited therein).
By the complementary-dual principle, all the fixed points  can be obtained analytically in terms of the canonical dual solutions.
Their stability and extremality are  identified by the triality theory.
Applications are illustrated by problems  governed by nonconvex polynomial, exponential and logarithmic functions.
Our examples   show that both globally stable and locally stable/unstable   fixed point problems in $\real^n$ can be   can be obtained
 easily  by solving the associated canonical dual problems in $\real^m$ with $m< n$.
 However, the local stability condition  for those fixed points $\barbx(\barbvsig)$ with indefinite  $\bG(\barbvsig)$  still remains unknown and it deserves seriously study in the future.
  Also, the results presented in this paper can be  generalized to  problems with nonsmooth potential functions.

\noindent{\bf Acknowledgement}:
The research was supported by US Air Force Office of Scientific Research
under the grants  (AOARD) FA2386-16-1-4082 and FA9550-17-1-0151.\\

\end{document}